\documentclass[12pt]{article}
\usepackage{bbm}
\usepackage{bbm}
\usepackage{amsfonts}
\usepackage{amsfonts}
\usepackage{amsfonts}
\usepackage{amsfonts}
\usepackage{amsfonts}
\usepackage{amsfonts}
\usepackage{bbm}
\usepackage{mathrsfs}
\usepackage{amssymb,amsfonts,amsmath,amsthm,colordvi,cite}
\usepackage{epsfig}
\newtheorem{thm}{Theorem}[section]

\newtheorem{lem}[thm]{Lemma}

\newtheorem{conj}[thm]{Conjecture}

\makeatletter \@addtoreset{equation}{section}
\def\qed{\hfill \rule{4pt}{7pt}}

\begin{document}

\begin{center}
{\bf\Large{2-Log-concavity of the Boros-Moll Polynomials  } }
\end{center}

\begin{center}
William Y. C. Chen$^1$ and Ernest X. W. Xia$^2$

$^1$Center for Combinatorics, LPMC-TJKLC\\
Nankai University\\
 Tianjin 300071, P. R. China

$^2$Department of Mathematics\\
 Jiangsu University\\
Jiangsu, Zhenjiang   212013, P. R. China

Email: $^1$chen@nankai.edu.cn, $^2$xxw@cfc.nankai.edu.cn

\end{center}

%======================================================

\noindent {\bf Abstract.} The Boros-Moll polynomials $P_m(a)$  arise
in the evaluation of a quartic  integral. It has been conjectured by Boros and Moll that these polynomials are infinitely log-concave. In this paper,
we show that $P_m(a)$ is $2$-log-concave for any $m\geq 2$.
 Let $d_i(m)$ be the coefficient of $a^i$ in $P_m(a)$.
We also show that the sequence  $\{ i
   (i+1)(d_i^{\,2}(m)-d_{i-1}(m)d_{i+1}(m)) \}_{1\leq i \leq m}$  is log-concave.
       This leads another proof
   of Moll's minimum conjecture.



\section{Introduction}

The objective of this paper  is to prove the
2-log-concavity  of  the Boros-Moll  polynomials.
 Recall that a sequence
 $\{a_i\}_{0 \leq i \leq  n}$
 of real  numbers is said to be
  unimodal  if
 there
exists an index $0
\leq j \leq n$ such that
\[
a_0\leq a_1\leq \cdots \leq a_{j-1}
 \leq a_j\geq a_{j+1}\geq \cdots
\geq a_n.
\]
Set $a_{-1}=0$ and $a_{n+1}=0$.
  We say that
    $\{a_i\}_{0 \leq i \leq  n}$
  is   log-concave if
\[
a_i^2-a_{i+1}a_{i-1}\geq 0,
 \qquad 1\leq i \leq n.
\]
  A polynomial is said to be  unimodal
  (resp., log-concave)
  if the sequence  of its
coefficients is unimodal (resp.,
 log-concave).
 It is easy to see
that for  a  positive sequence, the
 log-concavity is  stronger than  the unimodality.
For a sequence $A=\{a_i\}_{0 \leq i \leq n}$,
 we define the  operator $\mathcal {L}$
 by $\mathcal{L}(A)
 =\{b_i\}_{0\leq i \leq n}$,
  where
 \begin{align}\label{c1-operator}
b_i=a_i^2-a_{i-1}a_{i+1}, \quad  0\leq i \leq n.
 \end{align}
 We say that
$\{a_i\}_{0 \leq i \leq n}$ is
 $k$-log-concave  if the sequence
$\mathcal {L}^j \left(\{a_i\}_{0 \leq i \leq n}
 \right)$
 is log-concave
   for
every $0 \leq j \leq k-1$,
 and that $\{a_i\}_{0 \leq i \leq n}$ is
$\infty$-log-concave  if $\mathcal {L}^k \left(\{a_i\}_{0 \leq i
\leq n}\right)$ is
 log-concave
 for every $k \geq 0$.

Boros and Moll \cite{George2004}  conjectured that the binomial
coefficients ${n\choose k}$ are infinitely log-concave for any
$n$. An generalization of this conjecture was given
independently
by Fisk \cite{Fisk},
 McNamara
and  Sagan \cite{McNamara}, and Stanley, see \cite{Ba}, which
states that if a polynomial $a_0+a_1x+\cdots +a_nx^n$ has
only real zeros, then the polynomial $b_0+b_1x +\cdots +b_nx^n$ also
has only  real zeros, where $b_i=a_i^2-a_{i-1}a_{i+1}$.
This conjecture has been proved by
Br\"{a}nd\'{e}n \cite{Ba}. While Br\"{a}nd\'{e}n's theorem does
not directly apply to the Boros-Moll polynomials, the
$2$-log-concavity and $3$-log-concavity can be recasted in terms of
the real rootedness of certain  polynomials derived from the
Boros-Moll polynomials, as conjectured
 by Br\"{a}nd\'{e}n. It is worth mentioning that   McNamara
and  Sagan \cite{McNamara}
 conjectured that  for fixed $k$,
 the $q$-Gaussian coefficients ${n \brack k}$ are
   infinitely $q$-log-concave.  Chen,
  Wang  and Yang \cite{Chen-Wang-Yang} proved the strong $q$-log-concavity of  the $q$-Narayana numbers $N_q(n,k)$ for fixed $k$, which turns out to be
  equivalent to  the $2$-fold $q$-log-concavity of the Gaussian coefficients.

Recall that   Boros and Moll
 \cite{George1999-1,George1999-2,George1999-3,George2001
 ,George2004,Moll2002}
have studied the following
quartic integral and
  have shown that for any $a>-1$ and any
nonnegative integer $m$,
\begin{equation*}
 \int_{0}^\infty
\frac{1}{(x^4+2ax^2+1)^{m+1}}dx
=\frac{\pi}{2^{m+3/2}(a+1)^{m+1/2}}P_m(a),
\end{equation*}
 where
\begin{align}\label{E1}
P_m(a)=\sum_{j,k}{2m+1 \choose 2j}{m-j \choose k}{2k+2j \choose
k+j}\frac{(a+1)^j(a-1)^k}{2^{3(k+j)}}.
\end{align}
    Using  Ramanujan's Master Theorem, Boros and Moll
\cite{George2001,Moll2002} obtained the following formula for
$P_m(a)$:
\begin{align}\label{E2}
P_m(a)=2^{-2m}\sum_{k}2^k{2m-2k
 \choose m-k}{m+k \choose k}(a+1)^k,
\end{align}
  which implies that
  $P_m(a)$ is a polynomial in $a$ with positive coefficients.
 Chen, Pang and Qu
 \cite{Chen-Pang-Qu-2008}
 gave  a combinatorial argument to show that
 the double sum \eqref{E1} can be reduced to the single
  sum  \eqref{E2}.
 Let $d_i(m)$ be the coefficient of $a^i$ of $P_m(a)$, that is,
\begin{equation} \label{pma}
P_m(a)=\sum_{i=0}^md_i(m)a^i.
\end{equation}
For any $m$, $P_m(a)$ is  called
 a Boros-Moll polynomial, and the sequence
$\{d_i(m)\}_{0\leq i \leq m}$  is called a  Boros-Moll sequence.
 From \eqref{E2},  we know that $d_i(m)$ can be expressed as
\begin{align}\label{Defi}
d_i(m)=2^{-2m}\sum_{k=i}^m2^k{2m-2k
\choose m-k}{m+k \choose
k}{k\choose i}.
\end{align}
Many proofs of the above formula can be found in the
survey of Amdeberhan and   Moll \cite{Am}.

 \allowdisplaybreaks

Many combinatorial properties of $\{d_{i}(m)\}_{0 \leq i \leq m}$
have been studied.
 Boros and Moll \cite{George1999-2}
  proved that the sequence
$\{d_i(m)\}_{0 \leq i \leq m}$ is unimodal and the maximum element
appears in the middle. In other words,
\[d_0(m)<d_1(m)< \cdots<d_{\left[\frac{m}{2}\right]-1}(m)
<d_{\left[\frac{m}{2}\right]}(m)>
d_{\left[\frac{m}{2}\right]+1}(m)>\cdots >d_m(m).\] They
 also
   established the unimodality by   a
    different approach \cite{George1999-3}.
  Moll \cite{Moll2002}  conjectured that the sequence
  $\{d_i(m)\}_{0 \leq i \leq m}$ is
log-concave. Kauers and Paule \cite{Kauers2006} proved this
conjecture based on recurrence relations
 which were found by using  a  computer
algebra approach. Chen, Pang and Qu \cite{CPQ2010a} gave a
combinatorial
 proof of the log-concavity of $P_m(a)$
 by introducing the structure of partially $2$-colored  permutations.
  Chen and Gu \cite{Chen-Gu-2008} proved the
  reverse ultra log-concavity of the sequence $\{d_i(m)\}_{0 \leq i \leq m}$.
   Amdeberhan, Manna and Moll \cite{T-A}
studied the 2-adic valuation of an integer
sequence and found a combinatorial interpretation of the valuations
of the integer sequence which
 is related  to  the Boros-Moll sequences.
 Recently,  Chen and Xia  \cite{Chen2008}
 showed that the sequence  $\{d_i(m)\}_{0 \leq i \leq m}$
 satisfies the  strongly   ratio monotone
 property which implies the log-concavity and
 the spiral property. They  \cite{Chen-Xia-2009}
  also confirmed a conjecture
 of Moll which says that $\{i(i+1)\left(d_i^{\,2}(m)
 -d_{i-1}(m)d_{i+1}(m)\right)\}_{1\leq i \leq m}$
 attains  its minimum at $i=m$.

 Boros and Moll \cite{George2004} also
 made
  the following conjecture.

\begin{conj}\label{Moll-conjecture}
The Boros-Moll sequence $\{d_i(m)\}_{0\leq i \leq m}$ is
 $\infty$-log-concave.
\end{conj}

As shown by Boros and Moll  \cite{George1999-2}, in general,
$P_m(a)$ are not polynomials with only real zeros. Thus
 the theorem of Br\"{a}nd\'{e}n  \cite{Ba}  does  not   apply to $P_m(a)$.
Nevertheless,
  Br\"{a}nd\'{e}n  \cite{Ba} made the following
   conjectures on the real rootedness of polynomials
   derived from $P_m(a)$. These conjectures imply the 2-log-concavity and the 3-log-concavity of the
   Boros-Moll polynomials.

\begin{conj}[Br\"{a}nd\'{e}n] \label{P-conj-1}
For each positive integer $m$, the polynomial
\[
Q_m(x)=\sum_{i=0}^m \frac{d_i(m)}{i!}x^i
\]
has only real zeros.
\end{conj}

\begin{conj}[Br\"{a}nd\'{e}n]  \label{P-conj-2}
For each positive integer $m$, the polynomial
\[
R_m(x)=\sum_{i=0}^m \frac{d_i(m)}{(i+2)!}x^i
\]
has only real zeros.
\end{conj}

Note that $Q_m(x)=\frac{d}{dx^2} (x^2R_m(x))$. Hence $Q_m(x)$ has
only real zeros if $R_m(x)$ does.
 This yields that Conjecture \ref{P-conj-2} is stronger than
 Conjecture \ref{P-conj-1}. Based on a result of
 of Craven and
 Csordas \cite{Craven}, it can be seen that Conjecture \ref{P-conj-1}
 implies that  $P_m(a)$ is 2-log-concave and  Conjecture \ref{P-conj-2}
 implies that $P_m(a)$ is 3-log-concave.
   Conjectures \ref{P-conj-1} and \ref{P-conj-2}
     are still open.

In another direction, Kauers and Paule \cite{Kauers2006} considered
  using the approach
of recurrence relations to prove the $2$-log-concavity of $P_m(a)$, and
they indicated that there is little hope to make it work
 since the recurrence relations
are too complicated.

Roughly speaking, the main idea of this
 paper is to find an intermediate function $f(m,i)$
  so that we can reduce quartic inequalities
   for the $2$-log-concavity
   to quadratic inequalities.  To be precise,
the $2$-log-concavity is stated as follows.

\begin{thm}\label{Theo-1}
The Boros-Moll sequences are $2$-log-concave,
 that is, for $1\leq i\leq m-1$,
 \begin{align}\label{Th-1}
\frac{d_{i-1}^{\,2}(m)-d_{i-2}(m)d_{i}(m)}
{d_{i}^{\,2}(m)-d_{i-1}(m)d_{i+1}(m)} < \frac{d_{i}^{\,2}(m)
-d_{i-1}(m)d_{i+1}(m)} {d_{i+1}^{\,2}(m)-d_{i}(m)d_{i+2}(m)}.
 \end{align}
\end{thm}

The intermediate function $f(m,i)$ is given by
\begin{equation} \label{fmi}
f(m,i)=\frac{(i+1)(i+2)(m+i+3)^2}
{(m+1-i)(m+2-i)(m+i+2)^2}.
\end{equation}
Using this intermediate function, we can divide the
2-log-concavity into two quadratic inequalities, which are
stated below.

\begin{thm} \label{Theo-2} For
$1\leq i\leq m-1$, we have
\begin{align}\label{Th-2}
 \frac{(i+1)(i+2)(m+i+3)^2}
{(m+1-i)(m+2-i)(m+i+2)^2} <\frac{d_{i}^{\,2}(m) -d_{i-1}(m)d_{i+1}(m)}
{d_{i+1}^{\,2}(m)-d_{i}(m)d_{i+2}(m)}.
\end{align}
\end{thm}

\begin{thm} \label{Theo-3}
For $1\leq i\leq m-1$, we have
\begin{align}\label{Th-3}
\frac{d_{i-1}^{\,2}(m)-d_{i-2}(m)d_{i}(m)}
{d_{i}^{\,2}(m)-d_{i-1}(m)d_{i+1}(m)}
 < \frac{(i+1)(i+2)(m+i+3)^2}
{(m+1-i)(m+2-i)(m+i+2)^{2}}.
\end{align}
\end{thm}

As will be seen, the $2$-log-concavity of $P_m(a)$ implies
the  log-concavity
of a sequence considered by Moll \cite{Moll2007,Moll2008}.

\begin{thm} \label{Theo-4}
 For $m\geq 2$, the sequence
  $\{i(i+1) (d_i^{\,2}(m)-d_{i-1}(m)d_{i+1}(m)\}_{1\leq i \leq m}$ is
log-concave.
\end{thm}

Since log-concavity implies unimodality, the
above property leads to another proof of
Moll's minimum conjecture \cite{Moll2007} for the sequence
$\{i(i+1) (d_i^{\,2}(m)-d_{i-1}(m)d_{i+1}(m)\}_{ 1\leq i \leq m}$.
By comparing the first entry with the last
entry, we deduce  that
this sequence attains its minimum at $i=m$ which
 equals  $2^{-2m}m(m+1){2m \choose m}^2$.
 This conjecture was confirmed by  Chen and Xia \cite{Chen-Xia-2009}
 by using
 a result of Chen and Gu
\cite{Chen-Gu-2008}  and the spiral property of the Boros-Moll
 sequences
\cite{Chen2008}.

\section{How to guess the intermediate function $f(m,i)$}

In this section, we explain how we found the intermediate
function $f(m,i)$. We begin with
 a brief review of Kauers and Paule's approach
 to proving the log-concavity of  the Boros-Moll
polynomials \cite{Kauers2006}, because we need the
recurrence relations and an inequality established
by Kauers and Paule.  Here are the four recurrence relations
\begin{align}
d_i(m+1)=&\frac{m+i}{m+1}d_{i-1}(m)+\frac{(4m+2i+3)}{2(m+1)}d_i(m),
 \ \ \ \ 0 \leq i \leq m+1, \label{recu1}\\[6pt]
 d_{i}(m+1)=&\frac{(4m-2i+3)(m+i+1)}{2(m+1)(m+1-i)}d_i(m)
 \nonumber \\[6pt]
 & \qquad\qquad \qquad  -\frac{i(i+1)}{(m+1)(m+1-i)}d_{i+1}(m),
 \qquad \ \ \ 0 \leq i \leq
 m, \label{recu2}\\[6pt]
 d_i(m+2)=&\frac{-4i^2+8m^2+24m+19}{2(m+2-i)
 (m+2)}d_i(m+1) \nonumber \\[6pt]
 & \qquad\quad  -\frac{(m+i+1)(4m+3)(4m+5)}
 {4(m+2-i)(m+1)(m+2)}d_i(m), \qquad \ 0 \leq
  i \leq m+1,\label{recu3}
\end{align}
and for $0 \leq i \leq m+1$,
\begin{align}\label{recu4}
(m+2-i)(m+i-1)d_{i-2}(m)-(i-1)(2m+1)d_{i-1}(m)+i(i-1)d_i(m)=0.
\end{align}
These recurrences are
derived  by Kauers and Paule
\cite{Kauers2006}.
 In fact, the relations
 \eqref{recu3} and \eqref{recu4}
  are  derived independently by Moll
  \cite{Moll2007} via the WZ-method \cite{Wilf1992},
  and the other two relations
   (\ref{recu1}) and (\ref{recu2}) can
  be easily deduced
    from \eqref{recu3} and \eqref{recu4}.
    Based on the  four recurrence relations,
Kauers and Paule \cite{Kauers2006} proved the
following inequality
  from which the log-concavity of the Boros-Moll
  sequences can be deduced.

\begin{thm} \label{tkp}
{\rm (Kauers and Paule \cite{Kauers2006})} Let $m,i$ be integers
 with $m\geq 2$. For $0 < i< m$, we have
\begin{align}\label{condition}
\frac{d_i(m+1)}{d_i(m)} \geq \frac{4m^2+7m+i+3}{2(m+1-i)(m+1)}.
\end{align}
\end{thm}

Chen and Gu \cite{Chen-Gu-2008}
showed that
 $\{i!d_{i}(m)\}_{0 \leq i \leq
m}$ is log-concave and the sequence
 $\{d_{i}(m)\}_{0 \leq i \leq
m}$ is
 reverse ultra log-concave.
  They established  the following
 upper bound for $d_i(m+1)/d_i(m)$.

\begin{thm}  {\rm (Chen and Gu \cite{Chen-Gu-2008})}
\label{Upper-band}  Let $m,i$  be  integers and $m \geq 2$.
 We have
for
 $0 \leq i \leq m$,
\begin{align}\label{u-band}
\frac{d_i(m+1)}{d_i(m)}\leq \frac{4m^2+7m
+3+i\sqrt{4m+4i^2+1}-2i^2}{2(m+1)(m+1-i)}.
\end{align}
\end{thm}

 Theorems \ref{tkp} and \ref{Upper-band} are needed in the proofs
of Theorems
 \ref{Theo-2} and \ref{Theo-3}, and they are also needed to have
 a good guess of the intermediate function  $f(m,i)$.
  We start with an approximation of
 \[ \frac{d_{i-1}^{\,2}(m)-d_{i-2}(m)d_{i}(m)}
{d_{i}^{\,2}(m)-d_{i-1}(m)d_{i+1}(m)}.\]
  Recall that the following relation was proved by
 Chen and Gu \cite{Chen-Gu-2008},
\[
\lim_{m\rightarrow+\infty}
\frac{d_i^{\,2}(m) }{
\left(1+\frac{1}{i}\right) \left(1+\frac{1}{m-i}\right)
d_{i-1}(m)d_{i+1}(m)}=1.
\]
This implies that
\begin{equation}\label{di12m}
\frac{d_{i-1}^{\,2}(m)-d_{i-2}(m)d_{i}(m)}
{d_{i}^{\,2}(m)-d_{i-1}(m)d_{i+1}(m)} \approx
\frac{(i+1)(m+1-i)d_{i-1}^{\,2}(m)}
{i(m+2-i)d_i^{\,2}(m)}.
\end{equation}
Using the recurrence relation  \eqref{recu1}, we find
\begin{equation}\label{di-12}
\frac{d_{i-1}^{\,2}(m)} {d_i^{\,2}(m)} =\frac{(m+1)^2d_i^{\,2}(m+1)}
{(m+i)^2d_i^{\,2}(m)} - \frac{(4m+2i+3)(m+1)d_{i}(m+1)}
 {(m+i)^2d_i(m)}+\frac{(4m+2i+3)^2}{4(m+i)^2}.
\end{equation}
On the other hand, by Theorems  \ref{tkp} and \ref{Upper-band}, we
get
\begin{align*}
\lim_{m\rightarrow +\infty}\frac{2(m+1)(m+1-i)d_i(m+1)}
{(4m^2+7m+i+3)d_i(m)}=1.
\end{align*}
It follows that
\begin{equation}\label{dim}
\frac{d_i(m+1)}{d_i(m)}\approx \frac{4m^2+7m+i+3}
 {2(m+1)(m+1-i)}.
\end{equation}
Substituting (\ref{dim}) into (\ref{di-12}) yields
\begin{equation}\label{dim2}
{d_{i-1}^{\,2} (m) \over d_i^{\,2}(m) }
\approx\frac{i^2(i+1+m)^2}{(m+1-i)^2(m+i)^2}.
\end{equation}
Combining (\ref{di12m}) and (\ref{dim2}), we deduce that
\begin{equation} \label{try}
\frac{d_{i-1}^{\,2}(m)-d_{i-2}(m)d_{i}(m)}
{d_{i}^{\,2}(m)-d_{i-1}(m)d_{i+1}(m)}
 \approx \frac{i(i+1)(m+1+i)^2}
 {(m+1-i)(m+2-i)(m+i)^2}.
\end{equation}
It turns out that the above expression is not an intermediate function
that we are looking form. Naturally, we should try to make it a little
bigger. The above expression gives a guideline for a suitable adjustment.
Let us consider the shifts of
the factors in the expression (\ref{try}).
After a few trials,
we find that the  function below serves the purpose
as a desired intermediate function
\begin{equation}
\frac{(i+1)(i+2)(m+i+3)^2}
{(m+1-i)(m+2-i)(m+i+2)^2},
\end{equation}
which is the function $f(m,i)$ as given by
(\ref{fmi}).

\section{Proof of Theorem  \ref{Theo-2}}

In this section, we aim to
 give a proof of Theorem \ref{Theo-2}. The idea
goes as follows.
We wish to prove an equivalent form of Theorem \ref{Theo-2}, that is,
 the difference
\begin{align}
&(m+1-i)(m+2-i)(m+i+2)^2 \left(d_{i}^{\,2}(m)
-d_{i-1}(m)d_{i+1}(m)\right)\nonumber \\[6pt]
&\qquad\qquad\qquad
 -(i+1)(i+2)(m+i+3)^2\left(d_{i+1}^{\,2}(m)
-d_{i}(m)d_{i+2}(m)\right) \label{diff}
\end{align}
is positive.  As will be seen, in view of the recurrence relations of $d_i(m)$,
 (\ref{diff}) can be
written as
\begin{align}\label{form}
A(m,i)d_{i}^{\,^2}(m+1)+
B(m,i)d_{i}(m+1)d_i(m)+C(m,i)d_{i}^{\,^2}(m),
\end{align}
where $A(m,i)$, $B(m,i)$ and $C(m,i)$
are given by \eqref{A}, \eqref{B} and \eqref{C}.
To confirm that the quadratic form (\ref{form}) is
positive, we consider the quadratic polynomial in
$d_i(m+1)/d_i(m)$
\begin{align}\label{form-p}
A(m,i)\frac{d_{i}^{\,^2}(m+1)}{d_i^{\, 2}(m)}+
B(m,i)\frac{d_{i}(m+1)}{d_i(m)}+C(m,i).
\end{align}
It will be shown that
$A(m,i)<0$ for $1\leq i \leq m$. Moreover, we shall show
that the above polynomial has distinct real roots $x_1$ and
$x_2$. Assume that $x_1<x_2$. If the relation
\[ x_1 < \frac{d_{i}(m+1)}{d_i(m)} <x_2\]
holds, then  the quadratic polynomial (\ref{form-p})
is positive.

To present the following theorem, we need some notation.
Let
\begin{align}
A(m,i)&=-\frac{(m+1)^2(m+1-i)^2D(m,i)}
 {(m+i)i^2(i+1)},\label{A}\\[6pt]
 B(m,i)&=\frac{(i-m-1)(m+1)E(m,i)}
 {(i+m)i^2(i+1)},\label{B}\\[6pt]
 C(m,i)&=\frac{F(m,i)}{4(i+m)i^2(i+1)},\label{C}\\[6pt]
 \Delta_1(m,i)&=B^2(m,i) -4A(m,i)C(m,i)\nonumber \\[6pt]
 &=
\frac{(m+1-i)^2(m+1)^2\left(4(m+i)^2G(m,i)+
H(m,i)\right)}{i^2(i+m)^2(i+1)^2},\label{Delta-1}
\end{align}
where $D(m,i)$, $E(m,i)$, $F(m,i)$
 and $G(m,i)$ are given by
\begin{align*}
D(m,i)=&\;6m^2i+2m^2i^2
+21mi+14mi^2+4mi^3+10i\\[6pt]
&+17i^2+10i^3+2i^4+2m^3+12m^2+18m,\\[6pt]
E(m,i)=&\;4i^2(i^2-2m^2)(i+m)^2 + 2(i+m)(10i^4-4m^4-9im^3-27i^2m^2
-4i^3m)\\[6pt]
&+27i^4 -55i^3m -175i^2m^2 -139im^3
-62m^4-16i^3-155i^2m\\[6pt]
&-229im^2 -162m^3 -60i^2 -142im-162m^2-30i-54m,\\[6pt]
F(m,i)=&\;32i^2m^2(i-m)(i+m)^3+ 16m(4i^4+10i^3m-14i^2m^2-3im^3-2m^4)
(i+m)^2\\[6pt]
&+2(i+m)(-152m^5-250im^4-377i^2m^3
+111i^3m^2+181i^4m+15i^5) \\[6pt]
& +168i^5
 +694i^4m-280i^3m^2-2052i^2m^3
 -2160im^4 -1106m^5 +273i^4
\\[6pt]
&-i^3m-1809i^2m^2-2831im^3
-1968m^4+18i^3-898i^2m-1936im^2\\[6pt]
& -1836m^3-207i^2
 -663im-864m^2-90i-162m,\\[6pt]
 G(m,i)=&\;m^2(2i^3-m^2)^2
+(56i^6m-24i^3m^3)+(20i^5m^2-2i^2m^4)
 \\[6pt]
 &+4i^8+8i^7m
 +40i^7+169i^6+166i^5m
 +70i^4m^2,
\\[6pt]
H(m,i)=&\;1588i^7 +4440i^6m+4768i^5m^2+2148i^4m^3+324i^3m^4
+144i^2m^5
\\[6pt]
&+104im^6+52m^7
 +2345i^6+6666i^5m+6991i^4m^2+3624i^3m^3
 +1567i^2m^4\\[6pt]
& +646im^5+289m^6 +2418i^5 +7232i^4m+8044i^3m^2
+5340i^2m^3+2234im^4\\[6pt]
&+892m^5 +1903i^4+5810i^3m
+7225i^2m^2 +4104im^3+1618m^4
+1086i^3\\[6pt]
&+3332i^2m +3470im^2+1608m^3+321i^2+914im+657m^2.
\end{align*}

\begin{thm} \label{Th-root}
 For $1\leq i \leq m-1$ and $m\geq 126$, we have
\begin{align}
\frac{-B(m,i)
+\sqrt{\Delta_1(m,i)}}{2A(m,i)}<\frac{d_i(m+1)}{d_i(m)}
<\frac{-B(m,i)
-\sqrt{\Delta_1(m,i)}}{2A(m,i)}.\label{in-Th}
\end{align}
\end{thm}
In order to prove Theorem \ref{Th-root}, it is necessary to show that  $\Delta_1(m,i)>0$.

\begin{lem} \label{Lemma-3-1}
For $1 \leq i \leq m-1$ and $m\geq 126$,
  we have $\Delta_1(m,i)>0$.
\end{lem}

\noindent
{\it Proof.} In view of the definition \eqref{Delta-1}
 of $\Delta_1(m,i)$ and
 the fact that $H(m,i)$ is positive,
 it suffices to show that $G(m,i)>0$
 for $1 \leq i \leq m-1$. We consider three cases concerning the
 range of $i$. Case 1: $i^3 \geq \frac{3}{7}m^2$.
   In this case, we have
\[
m^2(2i^3-m^2)^2\geq 0,\quad
56i^6m-24i^3m^3\geq 0, \quad
20i^5m^2-2i^2m^4>0,
\]
and so $G(m,i)>0$.
Case 2: $\frac{m^2}{10}<
 i^3 <
\frac{3}{7}m^2$. In this case,
we have
\[
m^2(2i^3-m^2)^2\geq \frac{m^6}{49},\quad
 56i^6m-24i^3m^3\geq
-\frac{18}{7}m^5, \quad 20i^5m^2-2i^2m^4>0.
\]
Thus, for $m\geq 126$,
\[
G(m,i)\geq \frac{m^6}{49}-\frac{18}{7}m^5>0.
\]
Case 3: $1 \leq i^3\leq \frac{m^2}{10}$. In this case, we have
\[
m^2(2i^3-m^2)^2\geq \frac{16m^6}{25},\quad 56i^6m-24i^3m^3
 \geq -\frac{46}{25}m^5,\quad
20i^5m^2-2i^2m^4 >
 -2m^{16/3}.
\]
It follows that
\begin{align}\label{G}
G(m,i)\geq \frac{16m^6}{25}-\frac{46}{25}m^5
 -2m^{16/3}.
\end{align}
It is easily checked that the right-hand side of \eqref{G}
 is positive for $m\geq 10$. This completes the proof.\qed

We are now ready to prove Theorem \ref{Th-root}.

\noindent
{\it Proof of Theorem \ref{Th-root}.}
We first consider the lower bound of $d_{i}(m+1) / d_i(m)$, namely,
  \begin{align}\label{L-3-2}
\frac{d_{i}(m+1)}{d_i(m)}>\frac{-B(m,i)
 +\sqrt{\Delta_1(m,i)}}{2A(m,i)}.
\end{align}
From the inequality
  \eqref{condition} of Kauers and Paule \cite{Kauers2006}, we see that
  \eqref{L-3-2} is a consequence of the
  relation
\begin{align}\label{3-1}
\frac{4m^2+7m+i+3}{2(m+1)(m+1-i)}>\frac{-B(m,i)
+\sqrt{\Delta_1(m,i)}}{2A(m,i)}.
\end{align}
Since $A(m,i)<0$ for $1\leq i \leq m$, the inequality
 \eqref{3-1} can be rewritten as
\begin{align}\label{Left}
A(m,i)\frac{4m^2+7m+i+3}{(m+1)(m+1-i)}+B(m,i)
 <\sqrt{\Delta_1(m,i)}.
\end{align}
To verify \eqref{Left}, we calculate the difference of the squares
of both sides.
It is easily checked that
\begin{align*}
&\Delta_1(m,i)-\left(A(m,i) \frac{4m^2+7m+i+3}
{(m+1)(m+1-i)}+B(m,i)\right)^2 =\frac{(m+1-i)^2(m+1)^2K(m,i)}
{i^2(i+m)^2(i+1)^2},
\end{align*}
where $K(m,i)$ is given by
\begin{align*}
K(m,i)&=4(2i^4+4i^3m+2i^2m^2+10i^3+14i^2m+6im^2
+2m^3+17i^2+21im+12m^2\\[6pt]
&\qquad +10i+18m)(2i^3m^2 +2i^2m^3 -2i^5-2i^4m-9i^4
+2i^3m+16i^2m^2+6im^3\\[6pt]
&\qquad +m^4-7i^3+23i^2m+23im^2+9m^3
 +12i^2+16im +20m^2+8i+8m),
\end{align*}
which is positive for $1 \leq i \leq m-1$.
 Hence, by Lemma \ref{Lemma-3-1}, we obtain
  \eqref{Left}. This yields \eqref{L-3-2}.

It remains to  the consider the upper bound of $d_i(m+1)/d_i(m)$, namely,
\begin{align}\label{L-3-3}
\frac{d_{i}(m+1)}{d_i(m)}<\frac{-B(m,i)
 -\sqrt{\Delta_1(m,i)}}{2A(m,i)}.
\end{align}
By  Theorem  \ref{Upper-band} of Chen and
 Gu \cite{Chen-Gu-2008},
we see that (\ref{L-3-3}) is a consequence of the following relation
\begin{align}\label{3-2}
\frac{4m^2+7m+i\sqrt{4i^2+4m+1}-2i^2+3}
 {2(m+1)(m+1-i)}<\frac{-B(m,i)
 -\sqrt{\Delta_1(m,i)}}{2A(m,i)}.
\end{align}
Since $A(m,i)<0$  for $1\leq i \leq m-1$,
 \eqref{3-2} can be rewritten as
\begin{align}\label{Right}
A(m,i) \frac{4m^2+7m+i\sqrt{4i^2+4m+1}-2i^2+3}
 {(m+1)(m+1-i)}+B(m,i)>-\sqrt{\Delta_1(m,i)}.
\end{align}
As before, we can check  \eqref{Right} by
computing the difference of  the squares
of both sides. It is readily seen that
\begin{align*}
&\Delta_1(m,i)-\left(A(m,i) \frac{4m^2+7m+i\sqrt{4i^2+4m+1}-2i^2+3}
{(m+1)(m+1-i)}+B(m,i)\right)^2\\[6pt]
&\qquad\qquad
 =\frac{(m+1-i)^2(m+1)^2 L(m,i)}
  {i^2(i+m)(i+1)^2},
\end{align*}
where $L(m,i)$ is given by
\begin{align*}
L(m,i)=&\;2(2i^4+4i^3m+2i^2m^2+10i^3+14i^2m+6im^2
+2m^3+17i^2\\[6pt]
&\  +21im+12m^2+10i+18m)\big(-4i^3m-8i^2m^2
-4im^3 -20i^2m \\[6pt]
&\  -24im^2 -4m^3
+7i^2-28im-19m^2+20i-20m+7\\[6pt]
&\  +(2i^2m+4im^2+2m^3+i^2 +24m
+14im+13m^2+6i+9) \sqrt{4i^2+4m+1}\big).
\end{align*}
But
\begin{align*}
&(2i^2m+4im^2+2m^3+i^2+24m
 +14im+13m^2+6i+9)^2 (4i^2+4m+1)
 -(-4i^3m\\[6pt]
&\quad -8i^2m^2 -4im^3 -20i^2m
 -24im^2 -4m^3 +7i^2-28im-19m^2+20i-20m+7 )^2\\[6pt]
 =&16i^6m+96i^5m^2
 +176i^4m^3+128i^3m^4
 +48i^2m^5+32im^6+16m^7\\[6pt]
 &\quad+4i^6+264i^5m+972i^4m^2
 +1088i^3m^3+492i^2m^4+312im^5
 +196m^6\\[6pt]
 &\quad+48i^5+1456i^4m+3248i^3m^2
 +2064i^2m^3+1184im^4+960m^5
 +168i^4\\[6pt]
 &\quad+3508i^3m+4368i^2m^2
 +2372im^3+2384m^4+164i^3
 +3876i^2m\\[6pt]
 &\quad+3036im^2+3196m^3-120i^2
 +2164im+2404m^2-172i+1036m+32,
\end{align*}
which is positive for $1\leq i \leq m-1$. So we reach the conclusion that $L(m,i)>0$. Therefore, we obtain
  \eqref{Right} which implies \eqref{L-3-3}.
 In view of \eqref{L-3-2} and \eqref{L-3-3},
  we arrive at \eqref{in-Th}. This completes the proof. \qed

To conclude this section, we present a proof of Theorem \ref{Theo-2}.

\noindent
{\it Proof of Theorem \ref{Theo-2}.}
First we show that the difference
\eqref{diff} can be represented
in terms of $d_i(m)$ and $d_i(m+1)$. From
 \eqref{recu1},  \eqref{recu2} and  \eqref{recu4},
 it follows that for $1\leq i \leq m-1$,
\begin{align}
d_{i+1}(m)&=\frac{(4m-2i+3)(m+i+1)} {2i(i+1)}d_i(m)
-\frac{(m+1-i)(m+1)}{i(i+1)}d_i(m+1),
\label{L-1}\\[6pt]
d_{i+2}(m)&=\frac{2m+1}{i+2}d_{i+1}(m)
-\frac{(m-i)(m+i+1)}{(i+1)(i+2)}d_i(m),
\label{L-2}\\[6pt]
d_{i-1}(m)&=\frac{m+1}{m+i}d_i(m+1) -\frac{4m+2i+3}{2(m+i)}d_i(m).
\label{L-3}
\end{align}
Applying the above recurrence relations, we find
\begin{align}\label{represent-1}
&(m+1-i)(m+2-i)(m+i+2)^2 \left(d_{i}^{\,2}(m)
-d_{i-1}(m)d_{i+1}(m)\right)\nonumber\\[6pt]
&\qquad\qquad
 -(i+1)(i+2)(m+i+3)^2\left(d_{i+1}^{\,2}(m)
-d_{i}(m)d_{i+2}(m)\right)\nonumber\\[6pt]
& \qquad\qquad\qquad =A(m,i)d_{i}^{\,2}(m+1)
+B(m,i)d_i(m+1)d_i(m)+C(m,i)d_i^{\,2}(m).
\end{align}
It is easy to check  that
  Theorem \ref{Theo-2} holds for $2\leq m \leq 125$.
 By Theorem \ref{Th-root}, we conclude that
  the difference \eqref{form} is positive for $m\geq 126$ and $1\leq i \leq m-1$.  This
 completes the proof. \qed

\section{Proof of Theorem  \ref{Theo-3}}

This section is devoted to a proof of
Theorem \ref{Theo-3}.  The main steps can be described as follows.
To prove the theorem, we wish to show that the difference
 \begin{align}\label{diff-1}
&(i+1)(i+2)(m+i+3)^2
 \left(d_{i}^{\,2}(m)-d_{i-1}(m)d_{i+1}
 (m)\right)\nonumber \\[6pt]
&\qquad\qquad -(m+1-i)(m+2-i)(m+i+2)^2\left(
d_{i-1}^{\,2}(m)-d_{i-2}(m)d_{i}(m)\right)
\end{align}
is positive for $1\leq i \leq m-1$.
By the recurrence relations of $d_i(m)$,  the difference
 \eqref{diff-1} can be restated as
  \begin{align}\label{form-1}
U(m,i)d_{i}^{\,^2}(m+1)+
V(m,i)d_{i}(m+1)d_i(m)+W(m,i)d_{i}^{\,^2}(m),
\end{align}
where $U(m,i)$, $V(m,i)$ and $W(m,i)$
are given by \eqref{U}, \eqref{V} and \eqref{W}.
We need to consider five cases for the range of $i$.
The conclusion in each case implies that  \eqref{form-1}
is positive. Notice that the definition of $\Delta_2(m,i)$
is given in (\ref{delta-2}), which can be either positive
or negative depending on the range of $i$.

\noindent
Case 1:
 $1\leq i < \left(\frac{m^2}
 {2}\right)^{1/3}-m^{1/3} $. In this case, $\Delta_2(m,i)$ can be
 either nonnegative or negative. We need to consider the
 case when $\Delta_2(m,i)$ is nonnegative. Theorem \ref{thm-4-6} is
 established for this purpose.

 \noindent
  Case 2: $\left(\frac{m^2}
 {2}\right)^{1/3}-m^{1/3}\leq i
 \leq \left(\frac{m^2}
 {2}\right)^{1/3}$. In this case, we show that $\Delta_2(m,i)<0$.

  \noindent
  Case 3:  $\left(\frac{m^2}{2}\right)^{1/3}
 < i < m^{2/3}$. In this case, $\Delta_2(m,i)$ can be
 either nonnegative or negative. We establish Theorem
 \ref{Lemma-4-7} when $\Delta_2(m,i)$ is nonnegative.

 \noindent
 Case 4: $m^{2/3}\leq i
    \leq m-4$. We show that  $\Delta_2(m,i)> 0$ and give a new lower bound on the ratio
$d_i(m+1)/d_i(m)$
 which implies that \eqref{form-1} is positive.

 \noindent
 Case 5: $m-3\leq i\leq m-1$. We can directly
  verify that \eqref{form-1} is positive.

The following notation will be used in the statement of Theorem  \ref{thm-4-6}.
Let
\begin{align}
U(m,i)&=\frac{(m+1)^2(m+1-i)
 R(m,i)}{i(m+i)^2},\label{U}\\[6pt]
V(m,i)&=\frac{(m+1)S(m,i)}
 {i(m+i-1)(m+i)^2}, \label{V}\\[6pt]
W(m,i)&=\frac{T(m,i)}
 {4i(m+i-1)(m+i)^2}, \label{W}\\[6pt]
 \Delta_2(m,i)&=V^2(m,i)
 -4U(m,i)W(m,i)
=\frac{(m+1)^2X(m,i)} {i(m+i)^2(m+i-1)^2},\label{delta-2}
\end{align}
where $R(m,i)$,
 $S(m,i)$, $T(m,i)$
 and $X(m,i)$
 are given by
\begin{align*}
R(m,i)=&\;2i^2m^2+4mi^3+6im^2+14mi^2
+2i^4+10i^3\\[6pt]
&\qquad +21mi
+17i^2+2m^3+12m^2+18m+10i,\\[6pt]
S(m,i)=&\; 4i^2(i^2-2m^2)(i+m)^3 +2(8i^4-4i^3m-21i^2m^2
-9im^3-4m^4)(i+m)^2\\[6pt]
&\qquad +(i+m)(-54m^4-121im^3
-99i^2m^2-41i^3m+7i^4)
-41i^4\\[6pt]
&\qquad -98i^3m-187i^2m^2-262im^3-100m^4
-41i^3-51i^2m \\[6pt]
&\qquad -106im^2+25i^2 +45im+108m^2+30i+54m,\\[6pt]
T(m,i)=&\; 32i^2m^2(i+m)^4 + 16m(4i^4+18i^3m+18i^2m^2
+7im^3+2m^4)(i+m)^2\\[6pt]
&\qquad +2(i+m)(120m^5+414im^4+601i^2m^3
+523i^3m^2+199i^4m+15i^5)\\[6pt]
&\qquad +132i^5+850i^4m+1912i^3m^2 +2652i^2m^3 +2084im^4+562m^5+153i^4\\[6pt]
&\qquad +417i^3m+983i^2m^2+1307im^3+300m^4 -48i^3-328i^2m\\[6pt]
&\qquad -248im^2-432m^3 -177i^2-405im-540m^2-90i-162m,
\\[6pt]
X(m,i)=&\;16i^7m^4-16i^4m^6+4im^8
 +64i^8m^3 -24i^2m^7+16i^{11}+64i^{10}m
 +96i^9m^2
\\[6pt]
&\qquad +(128i^{10} +448i^9m+624i^8m^2
+448i^7m^3+160i^6m^4-100i^3m^6)
\\[6pt]
&\qquad +(372i^9 +1280i^8m+1868i^7m^2 +1256i^6m^3 +128i^5m^4 -240i^4m^5 )
\\[6pt]
&\qquad +(340i^8 +1712i^7m+2520i^6m^2
+620i^5m^3 -1132i^4m^4
 -1096i^3m^5
\\[6pt]
&\qquad -528i^2m^6)+
(3692i^2m-52im^7-16m^8-523i^7 -2i^6m-509i^5m^2
\\[6pt]
&\qquad -2584i^4m^3
-3749i^3m^4-2910i^2m^5 -635im^6-176m^7 -1416i^6
\\[6pt]
&\qquad -5048i^3m^3-5940i^2m^4
 -1810im^5-656m^6-586i^5-3890i^4m
  \\[6pt]
&\qquad -3588i^2m^3-667im^4
-688m^5+1240i^4+1054i^3m+2274i^2m^2
\\[6pt]
&\qquad +3216im^3+1104m^4 +1221i^3 +2896im^2+2160m^3
-3550i^5m\\[6pt]
&\qquad  -4508i^4m^2-268i^2-2525i^3m^2 +488im-432m^2-524i-1296m).
\end{align*}
Obviously, $U(m,i)$ is positive
 for $1\leq i \leq m-1$.

In  Case 1, we obtain the following
inequality.

\begin{thm}\label{thm-4-6}
If $\Delta_2(m,i) \geq 0$,
 we have for $1\leq i \leq
\left(\frac{m^2}{2}\right)^{1/3}-m^{1/3}$
 and $m\geq 15$,
\begin{align}\label{4-1}
\frac{d_i(m+1)}{d_i(m)}
<\frac{-V(m,i) -\sqrt{\Delta_2(m,i)}} {2U(m,i)}.
\end{align}
\end{thm}

\noindent
{\it Proof.} From the inequality
  \eqref{u-band} of Chen and Gu \cite{Chen-Gu-2008}, we see that \eqref{4-1}  can be deduced from the following relation
\begin{align}\label{L-4-6}
\frac{4m^2+7m+i\sqrt{4i^2+4m+1}+3-2i^2}
  {2(m+1)(m+1-i)}
<\frac{-V(m,i) -\sqrt{\Delta_2(m,i)}} {2U(m,i)}.
\end{align}
To prove  \eqref{L-4-6},  let
\begin{align*}
A_1(m,i)&=2(m+1)(m+1-i),\\[6pt]
B_1(m,i)&=4m^2+7m+3-2i^2,\\[6pt]
C_1(m,i)&=4i^2+4m+1.
\end{align*}
Clearly,  \eqref{L-4-6} can be restated as
\begin{align}\label{in-D}
D_1(m,i)> A_1(m,i) \sqrt{\Delta_2(m,i)}
+2iU(m,i)\sqrt{C_1(m,i)},
\end{align}
where $D_1(m,i)$ is given by
\begin{align*}
D_1(m,i)=&\;-V(m,i) A_1(m,i)-2U(m,i)
B_1(m,i)\\[6pt]
=&\frac{2(m+1)^2(m+1-i)(2m+1)(i^2-i+m+m^2)
(m+2+i)^2}{(i+m)^2(i+m-1)}.
\end{align*}
Hence $D_1(m,i)$ is positive for $1\leq i \leq m$.
Since  $D_1(m,i)$ is positive,
 the inequality \eqref{in-D} follows from the
  inequality
\begin{align}\label{in-D-1}
D_1^2(m,i)> \left( A_1(m,i) \sqrt{\Delta_2(m,i)}
+2iU(m,i)\sqrt{C_1(m,i)}\right)^2,
\end{align}
which can be rewritten as
\begin{align}\label{in-D-2}
E_1(m,i)>4iA_1(m,i)U(m,i)
\sqrt{\Delta_2(m,i)C_1(m,i)},
\end{align}
where $E_1(m,i)$
is given by
\begin{align}\label{E-1}
E_1(m,i)=D_1^2(m,i) -A_1^2(m,i) \Delta_2(m,i)- 4i^2U^2(m,i)
 C_1(m,i).
\end{align}
It can be seen that  \eqref{in-D-2} is valid if
 $E_1(m,i)$ is positive and  the following inequality holds,
 \begin{align}\label{in-D-3}
 E_1^2(m,i)> 16i^2A_1^2(m,i)U^2(m,i)
\Delta_2(m,i)C_1(m,i) .
 \end{align}
Given the definition \eqref{E-1} of $E_1(m,i)$,
  it is easily checked that
\begin{align}\label{4-2}
E_1(m,i)=-\frac{8(m+1-i)^2(m+1)^4 R_1(m,i) S_1(m,i)}
{i(m+i-1)(m+i)^3},
\end{align}
where $R_1(m,i)$ and $S_1(m,i)$
 are given by
\begin{align}
R_1(m,i)=&2i^2m^2
+4mi^3+6im^2+14mi^2+2i^4 \nonumber\\[6pt]
&\qquad +10i^3 +21mi+17i^2+2m^3+12m^2+18m+10i, \label{R-1}\\[6pt]
S_1(m,i)=&8i^5m^2-4i^2m^4+36i^4m^2
 +12i^3m^3+(16i^6m-4m^5)+(8i^7-2im^4)\nonumber\\[6pt]
 &\qquad +(32i^6+52i^5m +30i^5+88i^4m
+66i^3m^2 -28i^2m^3-6im^4)\nonumber\\[6pt]
&\qquad +(36m -27i^4+55i^3m-65i^2m^2-23im^3-24m^4-56i^3
\nonumber\\[6pt]
&\qquad -101i^2m-9im^2-32m^3-9i^2 -20im+24m^2+22i).\label{S-1}
\end{align}
Using the expression \eqref{4-1} of $ E_1(m,i)$,
we see that the positivity of $ E_1(m,i)$ can be derived
 from the fact  that $S_1(m,i)$ is negative
  for $1\leq i \leq
\left(\frac{m^2}{2}\right)^{1/3}-m^{1/3}$
 and $m\geq 15$. We now proceed to show that $S_1(m,i)$ is
   negative. For $15 \leq m\leq 728$,
    the  claim  can be directly verified. Therefore,
     we may assume
    that $m\geq 729$.
  By putting the terms of $S_1(m,i)$ into groups as
  given in \eqref{S-1}, it is straightforward to see that
   the sum in every pair of parentheses  in \eqref{S-1} is negative for  $1\leq i \leq
  \left(\frac{m^2}{2}\right)^{1/3} -m^{1/3}$ and
  $m\geq 729$.
   Moreover, we can check that
\[
8i^5m^2-4i^2m^4<-15m^{11/3}i^2
+20m^{10/3}i^2-8m^3i^2.
\]It follows that
\begin{align*}
S_1(m,i)&<-15m^{11/3}i^2 +20m^{10/3}i^2 -8m^3i^2+36i^4m^2
 +12i^3m^3\\[6pt]
 &<(-5m^{5/3}+43m^{4/3})m^2i^2,
\end{align*}
which is negative when $m\geq 729$.
So we conclude that  $ E_1(m,i)>0$ for $1\leq i \leq \left(\frac{m^2}{2}
\right)^{1/3} -m^{1/3}$ and $m
\geq 15$.

We now turn to the proof of \eqref{in-D-3}.
Consider the  difference of the squares of both sides.
 It is  routine to check that
\begin{align}\label{F-1}
F_1(m,i)&=E_1^2(m,i)
 -16i^2 U^2(m,i)A_1^2(m,i)
\Delta_2(m,i)C_1(m,i)\nonumber\\[6pt]
&=\frac{-256(m+1-i)^4(m+1)^8 M_1^2(m,i)N_1(m,i)} {
i^2(i+m-1)^2 (i+m)^6 },
\end{align}
where $M_1(m,i)$ and $N_1(m,i)$
 are given by
 \begin{align*}
M_1(m,i)= &2i^4+4i^3m+2i^2m^2+10i^3 +14i^2m
+6im^2\\[6pt]
&  +2m^3+17i^2
+21im+12m^2+10i+18m,\\[6pt]
N_1(m,i)= &4i^{10}m-40i^8m^3-96i^7m^4 -128i^6m^5-128i^5m^6-88i^4m^7
-32i^3m^8-4i^2m^9\\[6pt]
&+i^{10}+12i^9m-92i^8m^2-400i^7m^3
-774i^6m^4-1100i^5m^5-1072i^4m^6\\[6pt]
&-592i^3m^7 -171i^2m^8-32im^9-4m^{10}+6i^9-58i^8m
-556i^7m^2-1602i^6m^3\\[6pt]
&-3236i^5m^4-4334i^4m^5-3204i^3m^6-1270i^2m^7
-322im^8-48m^9-3i^8\\[6pt]
&-351i^7m-1487i^6m^2 -4194i^5m^3-7663i^4m^4-7213i^3m^5-3519i^2m^6
\\[6pt]
&-1122im^7-208m^8-87i^7-695i^6m-2422i^5m^2
-5984i^4m^3-6495i^3m^4\\[6pt]
&-3165i^2m^5-1272im^6-336m^7-161i^6
-399i^5m-1212i^4m^2-107i^3m^3\\[6pt]
&+2447i^2m^4+1012im^5+104m^6+87i^5+839i^4m
+3175i^3m^2+6101i^2m^3\\[6pt]
&+2902im^4+816m^5+377i^4+1388i^3m+3137i^2m^2
+862im^3+432m^4\\[6pt]
&+32i^3-20i^2m-1308im^2 -432m^3-252i^2
 -720im-324m^2.
 \end{align*}
It is now easy to see that $N_1(m,i)<0$
 for $1 \leq i  < \left(\frac{m^2}{2}
  \right)^{1/3}-m^{1/3}$ and $m\geq 15$.
 So we have $F_1(m,i)>0$  for $1 \leq i  < \left(\frac{m^2}{2}
  \right)^{1/3}-m^{1/3}$ and $m\geq 15$. Hence the inequality \eqref{in-D-3} holds.  This completes the proof. \qed

For Case 2, the following lemma asserts that $\Delta_2(m,i)$ is negative.

\begin{lem}\label{Lem41}
For $ \left(\frac{m^2}{2}\right)^{1/3}-m^{1/3}
 \leq i
 \leq
 \left(\frac{m^2}{2}\right)^{1/3}$ and
 $m\geq 50$,
  we have $\Delta_2(m,i)<0$.
\end{lem}
\noindent{\it Proof.} By the definition \eqref{delta-2} of $\Delta_2(m,i)$,   it suffices to show that
$X(m,i)$
 is negative for $
\left(\frac{m^2}{2}\right)^{1/3}-m^{1/3} \leq i
 \leq
 \left(\frac{m^2}{2}\right)^{1/3}$
 and $m\geq 50$. For $50 \leq m \leq 2743$, the lemma  can be
  directly verified.
  Hence we may assume that $m\geq 2744$.
   Note that the expression in
    every pair of parentheses  is
negative for $\left(\frac{m^2}{2}
\right)^{1/3}-m^{1/3} \leq i
 \leq
 \left(\frac{m^2}{2}\right)^{1/3}$
  and $m\geq 2744$.
On the other hand, it can be checked that
\begin{align*}
&16i^7m^4-16i^4m^6+4im^8
= 4im^4(2i^3-m^2)^2 <58im^{22/3} \leq
47m^8,\\[6pt]
&64i^8m^3-24i^2m^7 +16i^{11}+64i^{10}m
 +96i^9m^2\leq -8i^2m^7+176i^9m^2
\leq -5m^{25/3} +22m^8.
\end{align*}
This yields
 \[
X(m,i)<-5m^{25/3} +69m^8.
 \]
 But the right-hand side of the above inequality is
 negative when $m\geq 2744$. This completes the proof.
\qed

As will be seen, Theorems 4.2 and 4.3 have the same
 expression of the lower bound for $d_i(m+1)/d_i(m)$.
 This expression will be needed in the proof of Theorem 1.6.
 It should be noted that for the case of Theorem 4.2,
 we shall show that this lower bound  can be derived
 from the lower bound of Kauers and Paule \cite{Kauers2006}.
 Numerical evidence shows that the bound in
 Theorem 4.3 seems sharper than the bound of Kauers and Paule when
 $i$ is large. However, we shall not make a rigorous comparison of these two bounds.

 For Case 3, we  have the following inequality.
 It should remarked that in this case $\Delta_2(m,i)$
 can be either positive or negative, and there is no
 need to specify the range of $i$
  for which $\Delta_2(m,i)$ is positive.

\begin{thm}\label{Lemma-4-7}
If $\Delta_2(m,i)\geq 0$, we have for  $\left(\frac{m^2} {2}\right)^{1/3} \leq i \leq
m^{2/3}$ and $m\geq 2$,
\begin{align}\label{L-4-7}
\frac{d_i(m+1)}{d_i(m)}
>\frac{-V(m,i)+\sqrt{\Delta_2(m,i)}}
{2U(m,i)}.
\end{align}
\end{thm}

\noindent
{\it Proof.} By the lower bound of $d_i(m+1)/d_i(m)$,
as given in  \eqref{condition}, we see that \eqref{L-4-7}  can be obtained from the following relation
\begin{align}\label{Kauers}
 \frac{4m^2+7m+i+3} {2(m+1)(m+1-i)}>\frac{-V(m,i)+\sqrt{\Delta_2(m,i)}}
{2U(m,i)},
\end{align}
which can be rewritten as
\begin{align}\label{in-E}
U(m,i) \frac{4m^2+7m+i+3}
 {(m+1)(m+1-i)} +
V(m,i)>\sqrt{\Delta_2(m,i)}.
\end{align}
In order to prove \eqref{in-E}, we shall show that
 for $\left(\frac{m^2}
 {2}\right)^{1/3} \leq i \leq
m^{2/3}$ and $m\geq 2$,
\begin{align}\label{I-1}
U(m,i) \frac{4m^2+7m+i+3}
 {(m+1)(m+1-i)} +
V(m,i)>0.
\end{align}
and
\begin{align}\label{I-2}
\left(U(m,i) \frac{4m^2+7m+i+3}
 {(m+1)(m+1-i)} +
V(m,i)\right)^2-\Delta_2(m,i)>0
\end{align}
We first deal with   inequality \eqref{I-1}.
 It is easily checked that
\begin{align*}
U(m,i) \frac{4m^2+7m+i+3}
 {(m+1)(m+1-i)} +
V(m,i)
 =\frac{(m+1)P(m,i)} {(m+i)^2(m+i-1)},
\end{align*}
where $P(m,i)$ is given by
\begin{align*}
P(m,i)=&\;4i^6 +(4i^3m^3-2m^5)
+(38i^3m^2-9m^4)+(14i^2m^3-11m^3)+12i^5m
\nonumber\\[6pt]
&\qquad +18i^5+44i^4m
 +(21i^4-10i^3)+60i^3m
+(35i^2m-21im)+12i^4m^2\nonumber\\[6pt]
& \qquad +(64i^2m^2-10m^2-22m)+16im^3 +(34im^2 -27i^2-6i).
\end{align*}
Since the sum in  every pair of parentheses
 in  the above  expression of $P(m,i)$
  is nonnegative  for $\left(\frac{m^2}
 {2}\right)^{1/3} \leq i \leq
m^{2/3}$ and $m\geq 2$,  it follows that $P(m,i)>0$. Thus,
 we obtain  \eqref{I-1}.

  We still need to consider  the inequality \eqref{I-2}.
 Clearly,
\begin{align*}
\left(U(m,i) \frac{4m^2+7m+i+3}
 {(m+1)(m+1-i)} +
V(m,i)\right)^2-\Delta_2(m,i)
 =\frac{4(m+1)^2
G_1(m,i)H_1(m,i)}
 {(m+i)^4(i+m-1)i},
\end{align*}
where $G_1(m,i)$ and $H_1(m,i)$ are given by
\begin{align*}
G_1(m,i)=&\;2i^4+4i^3m +2i^2m^2 +10i^3 +14i^2m
+6im^2\\[6pt]
&\qquad +2m^3+17i^2
+21im+12m^2+10i+18m,\\[6pt]
H_1(m,i)=&\;2i^7+4i^6m +7i^6 +11i^5m+8i^4m^2
+14i^3m^3+15i^2m^4+3i^4\\[6pt]
&\qquad +(7im^5-4i^4m^3)+(2m^6-2i^3m^4)
+7i^5+34i^4m+68i^3m^2+58i^2m^3
\\[6pt]
&\qquad +(29im^4-10i^2m)+(12m^5-12m^3)
 +(61i^3m-14i^2-40im)
\\[6pt]
& \qquad +(63i^2m^2-25im^2-18m^2) +21im^3+16m^4-5i^3.
\end{align*}
We see that  $G_1(m,i)>0$
  and
$H_1(m,i)>0$ for $\left(\frac{m^2}
 {2}\right)^{1/3} \leq i \leq
m^{2/3}$  and $m\geq 2$. Hence the  inequality \eqref{I-2}  holds.
 This completes the proof.  \qed

For Case 4, we give a lower bound for $d_i(m+1)/d_i(m)$ that takes
the same form as the lower bound in Case 3.

\begin{thm}\label{thm-4-4}
For $m\geq 273$ and $m^{2/3}\leq i \leq m-4$, we have
\begin{align}\label{L-4-4}
\frac{d_i(m+1)} {d_i(m)}> \frac{-V(m,i)
+\sqrt{\Delta_2(m,i)}}{2U(m,i)}.
\end{align}
\end{thm}

 For the clarity of presentation,   we   establish two lemmas for
 the proof of Theorem \ref{thm-4-4}.  First, we prove
 the  positivity
 of $\Delta_2(m,i)$ .

\begin{lem}\label{Lemma-4-1}
For $m^{2/3}\leq i \leq m-1$ and $m\geq 19$, we have
$\Delta_2(m,i)>0$.
\end{lem}
\noindent
{\it Proof.} By the definition \eqref{delta-2} of $\Delta_2(m,i)$, it suffices
 to show that $ X(m,i)$  is positive for
$m^{2/3}\leq i \leq m-1$ and $m\geq 19$.
By direct computation we find that the lemma holds
for $19\leq m \leq 132$.
 Moreover,  for $m\geq 133$ and
$m^{2/3}\leq i \leq m-1$, we have
\begin{align*}
X(m,i)&\geq 16i^{11}+64i^{10}m +96i^9m^2
+40i^8m^3+24(i^8m^3-i^2m^7)+16(i^7m^4-i^4m^6)\\[6pt]
&\quad +128i^{10}+(448i^9m-176m^7) +624i^8m^2
+(292i^7m^3-240i^4m^5-52im^7)
\\[6pt]
&\quad +(116i^6m^4-100i^3m^6-16m^8)
 +(1868i^7m^2 -1132i^4m^4-635im^6)
 \\[6pt]
&\quad +1096(i^6m^3-i^3m^5)+(160m^{2/3}-2910)i^2m^5
 +(620m^{2/3}-2584)i^5m^3
 \\[6pt]
 & \quad   +(128m^{4/3}-3749)i^3m^4+(4m-528)m^6i^2
 +(340m^{2/3}-523)i^7+1712i^7m
\\[6pt]
& \quad +(2520i^6m^2-2i^6m-509i^5m^2) +372i^9 +1280i^8m
-22928m^6-11944m^5\\[6pt]
&\geq 96m^8-22928m^6-11944m^5,
\end{align*}
which  is positive for $m\geq 133$.
This completes the proof. \qed

The proof of Theorem \ref{thm-4-4} is by induction on $m$.
The inductive argument requires an inequality concerning the
desired lower bound. We present this inequality in Lemma
\ref{Lemma-4-3},
and we need the following notation. Let
\begin{align*}
Y_1(m,i)&=\;\frac{(m+i+1)(4m+3)(4m+5)}
{4(m+2-i)(m+1)(m+2)},\\[6pt]
Y_2(m,i)& =\;\frac{-4i^2+8m^2+24m+19}
 {2(m+2-i)(m+2)},\\[6pt]
Y_3(m,i)&=\;2U(m+1,i)Y_2(m,i)
 +V(m+1,i)=\frac{(m+2)Y_5(m,i)}
 {(m+i)i(m+i+1)},\\[6pt]
Y_4(m,i)&=Y_3^2(m,i) -\Delta_2(m+1,i) =\frac{(m+2)^2Y_6(m,i)}
{(m+1+i)^2i^2(m+i)},
\end{align*}
where $Y_5(m,i)$ and $Y_6(m,i)$
  are  given by
\begin{align*}
Y_5(m,i)=&\;4i^2(2m^2-i^2)(i+m)^2 +2(i+m) (4m^4+7im^3
+31i^2m^2+4i^3m-12i^4)\\[6pt]
&\  -35i^4+59i^3m+199i^2m^2+151im^3
+82m^4+16i^3+181i^2m+321im^2\\[6pt]
&\  +282m^3 +70i^2+294im+368m^2+106i+160m,\\[6pt]
Y_6(m,i)=&\; (2i^4+4i^3m+2i^2m^2+14i^3+18i^2m+6im^2+2m^3
+33i^2+33im+18m^2\\[6pt]
&\  +37i+48m+32)\big(32i^2m^2(m-i)(i+m)^2 +16m(i+m)
 (2m^4
+im^3\\[6pt]
&\  +16i^2m^2-11i^3m-4i^4) -30i^5-394i^4m
-110i^3m^2+762i^2m^3+300im^4\\[6pt]
&\  +368m^5-168i^4-338i^3m+1154i^2m^2+558im^3
+1538m^4+1028i^2m\\[6pt]
&\  -141i^3 +631im^2+2882m^3 +391i^2+639im +2480m^2+260i+800m\big).
\end{align*}
 It is easily seen that  $Y_1(m,i)$, $Y_2(m,i)$,
 $Y_3(m,i)$ and $Y_4(m,i)$ are all
 positive
for $1\leq i \leq m-1$ and $m\geq 2$.

\begin{lem} \label{Lemma-4-3}
For $m^{2/3}\leq i \leq m-4$ and $m\geq 273$, we have
\begin{align}\label{L-4-3}
\frac{-V(m,i)+\sqrt{\Delta_2(m,i)}} {2U(m,i)}
>\frac{Y_1(m,i)}{Y_2(m,i)
-\frac{-V(m+1,i)+\sqrt{\Delta_2(m+1,i)}} {2U(m+1,i)}}.
\end{align}
\end{lem}

\noindent
{\it Proof.} Let us rewrite \eqref{L-4-3}  as
\begin{align}\label{in-F-1}
\frac{-V(m,i)+\sqrt{\Delta_2(m,i)}}
 {2U(m,i)}>
 \frac{2U(m+1,i)Y_1(m,i)}
 {Y_3(m,i)
 -\sqrt{\Delta_2(m+1,i)}}.
\end{align}
Since $Y_3(m,i)>0$ and $Y_4(m,i)>0$
 for $m^{2/3}\leq i \leq m-4$ and $m\geq 273$,  the inequality \eqref{in-F-1} follows from the inequality
   \begin{align}\label{in-F-2}
V(m,i)\sqrt{\Delta_2(m+1,i)} +Y_3\sqrt{\Delta_2(m,i)}
>Z_1(m,i)+\sqrt{\Delta_2(m,i)\Delta_2(m+1,i)},
 \end{align}
 where $Z_1(m,i)$ is given by
 \begin{align}
 Z_1(m,i)=&\;4U(m,i) U(m+1,i)
 Y_1(m,i)
 +V(m,i)Y_3(m,i). \label{Z-1}
 \end{align}
 Clearly, $Z_1(m,i)<0$ for
 $m^{2/3}\leq i \leq m-4$ and $m\geq 273$.
  To confirm  \eqref{in-F-2}, we shall show that
   the following three inequalities hold,
\begin{align}
&Z_1(m,i) +\sqrt{\Delta_2(m,i)\Delta_2(m+1,i)}<0, \label{ineq-2}\\[6pt]
&V(m,i)\sqrt{\Delta_2(m+1,i)} +Y_3(m,i)\sqrt{\Delta_2(m,i)}<0 \label{ineq-3}
\end{align}
and
\begin{align}
&\left(V(m,i)\sqrt{\Delta_2(m+1,i)}
+Y_3(m,i)\sqrt{\Delta_2(m,i)}\right)^2\nonumber\\[6pt]
&\qquad\qquad \qquad \qquad <\left(Z_1(m,i)
+\sqrt{\Delta_2(m,i)\Delta_2(m+1,i)}\right)^2.\label{ineq-1}
\end{align}
We first consider   inequality \eqref{ineq-2}. Let
\begin{align}
Z_2(m,i)=\Delta_2(m,i)\Delta_2(m+1,i)
 -Z_1^2(m,i).\label{Z-2}
\end{align}
Employing the  same  argument as in the proofs of  Lemmas \ref{Lemma-3-1}, \ref{Lem41}  and \ref{Lemma-4-1}, we find that $Z_2(m,i)<0$
 for $m^{2/3}\leq i \leq m-4$ and $m\geq 273$. The detailed proof is
 omitted since the expansion of $Z_2(m,i)$
 occupies more than three pages. Thus we obtain \eqref{ineq-2}
 since both $Z_1(m,i)$ and
   $Z_2(m,i)$ are negative for
 $m^{2/3}\leq i \leq m-4$ and $m\geq 273$.

We now turn to the proof of \eqref{ineq-3}. Note that $V(m,i)<0$ for $1\leq i \leq m-1$. Let
\begin{align}
Z_3(m,i)=Y_3^2(m,i)
 \Delta_2(m,i)
-V^2(m,i)\Delta_2(m+1,i).
\label{Z-3}
\end{align}
It is not difficult to show that $Z_3(m,i)<0$
 for
 $m^{2/3}\leq i \leq m-4$ and $m\geq 273$. But the detailed proof is
 omitted since the expansion of $Z_3(m,i)$ is too long.
   Since $Z_3(m,i)$,
   $V(m,i)$ are all negative  and $Y_3(m,i)$, $\Delta_2(m,i)$
    are  positive for
 $m^{2/3}\leq i \leq m-4$ and $m\geq 273$, we arrive at
  \eqref{ineq-3}.

It remains to prove   \eqref{ineq-1}, which can be restated as
\begin{align}\label{4-3}
Z_4(m,i)>Z_5(m,i)\sqrt{\Delta_2(m,i)\Delta_2(m+1,i)},
\end{align}
where $Z_4(m,i)$ and $Z_5(m,i)$ are given by
\begin{align}
  Z_4(m,i)=&\;
  V^2(m,i)\Delta_2(m+1,i)
  +Y_3^2(m,i)\Delta_2(m,i) \nonumber
  \\[6pt]
 &\qquad\qquad -Z_1^2(m,i)
 -\Delta_2(m,i)\Delta_2(m+1,i),\label{Z-4}\\[6pt]
 Z_5(m,i)=&\; 2Z_1(m,i) -2V(m,i) Y_3(m,i).\label{Z-5}
\end{align}
Using the same argument as in the proofs of Lemmas \ref{Lemma-3-1},
 \ref{Lem41}  and \ref{Lemma-4-1}, we can deduce that $Z_4(m,i)$ and $Z_5(m,i)$
  are   positive for
 $m^{2/3}\leq i \leq m-4$ and $m\geq 273$. Therefore,  \eqref{4-3}
  is a consequence of the fact that
  \begin{align}
Z_6(m,i) =\;Z_5^2(m,i)\Delta_2(m,i)\Delta_2(m+1,i)
 -Z_4^2(m,i)\label{Z-6}
\end{align}
is positive for
 $m^{2/3}\leq i \leq m-4$ and $m\geq 273$,
 which is not difficult to prove although $Z_6(m,i)$ is rather
  tedious. This completes the proof.
 \qed

We are now in a position to prove Theorem \ref{thm-4-4}.

\noindent
{\it Proof of Theorem \ref{thm-4-4}.}
 We proceed by induction on $m$. It is easy to
check that the theorem holds for $m=273$.
 We assume that the theorem
 is true for $n \geq 273$,  that is,
 \begin{align}\label{Assume}
d_{i}(n+1)\geq \frac{-V(n,i) +\sqrt{\Delta_2(n,i)}}{2U(n,i)}d_i(n),
 \qquad n^{2/3} \  \leq i \leq n-4.
 \end{align}
We aim to show that   \eqref{L-4-4} holds for $m=n+1$, that is,
 \begin{align}\label{Aim}
d_{i}(n+2)\geq \frac{-V(n+1,i) +\sqrt{\Delta_2(n+1,i)}}{2U(n+1,i)}
d_i(n+1) , \quad
 (n+1)^{2/3} \leq i \leq n-3.
 \end{align}
In view of Lemma \ref{Lemma-4-3} and
  inequality
  \eqref{Assume},  we find
\begin{align*}
d_i(n+1)>\frac{Y_1(n,i)} {Y_2(n,i) -\frac{-V(n+1,i)
+\sqrt{\Delta_2(n+1,i)}} {2U(n+1,i)}}d_i(n).
\end{align*}
It follows that for $n^{2/3}\leq i \leq n-4$,
\begin{align}\label{D-1}
Y_2(n,i)d_i(n+1)-Y_1(n,i)d_i(n)
>\frac{-V(n+1,i)+\sqrt{\Delta_2(n+1,i)}}
{2U(n+1,i)}d_i(n+1).
\end{align}
By the recurrence  relation \eqref{recu3},
the left hand side of
\eqref{D-1} equals  $d_i(n+2)$.
Thus we have verified  \eqref{Aim}  for
 $  (n+1)^{2/3}\leq i
\leq n-4$. It is still necessary to
show that \eqref{Aim} is true
for $i=n-3$, that is,
\begin{align}\label{Aim-2}
d_{n-3}(n+2)> \frac{-V(n+1,n-3) +\sqrt{\Delta_2(n+1,n-3)}}
{2U(n+1,n-3)} d_{n-3}(n+1).
\end{align}
Let
\begin{align*}
f(n)=&\;256n^{11}-4608n^{10}
+36544n^9-177920n^8+572592n^7-1218432n^6
\\[6pt]
&+1573768n^5-940352n^4 -66903n^3
-65525n^2-3657n-963.
\end{align*}
By the expression  \eqref{Defi} of $d_i(m)$,
we have
\begin{align*}
\frac{d_{n-3}(n+2)}{d_{n-3}(n+1)} =
&\;\frac{(2n+5)(16n^4+80n^3
+180n^2+240n+189)(2n-1)}
 {10(n+2)(45+72n+68n^2+48n^3+16n^4)}\\[6pt]
 >&\frac{12-65n+14n^2+3108n^4-3041n^3
 -1020n^5+136n^6+16n^7}{10(n+2)
 (2n-3)(1+2n+33n^2+4n^4-16n^3)}\\[6pt]
&\qquad +\frac{(n-1)\sqrt{(n-3)f(n)}}
{10(n+2)(2n-3)(1+2n+33n^2+4n^4-16n^3)}\\[6pt]
=&\frac{-V(n+1,n-3) +\sqrt{\Delta_2(n+1,n-3)}} {2U(n+1,n-3)}.
\end{align*}
Hence the proof is complete by induction. \qed

Finally, we are ready to complete the proof of Theorem \ref{Theo-3}.

\noindent
{\it Proof of Theorem \ref{Theo-3}.}
 For $2 \leq m \leq 272$,  the theorem
  can be easily verified. So we may assume that $m\geq 273$.
   The
difference \eqref{diff-1} can be represented
in terms of $d_i(m+1)$ and $d_i(m)$. From
 \eqref{recu4}  it follows
that
\begin{align}\label{L-4}
d_{i-2}(m)=\frac{(i-1)(2m+1)}{(m+2-i)(m+i-1)} d_{i-1}(m)
-\frac{i(i-1)}{(m+2-i)(m+i-1)}d_i(m).
\end{align}
Using recurrence relations \eqref{L-1}, \eqref{L-3}
 and  \eqref{L-4}, we find that
\begin{align}\label{represent-2}
&(i+1)(i+2)(m+i+3)^2
 \left(d_{i}^{\,2}(m)-d_{i-1}(m)d_{i+1}
 (m)\right)\nonumber\\[6pt]
&\qquad-(m+1-i)(m+2-i)(m+i+2)^2\left(
d_{i-1}^{\,2}(m)-d_{i-2}(m)d_{i}(m)\right)
\nonumber\\[6pt]
&\qquad\qquad =U(m,i)d_{i}^{\,2}(m+1) +V(m,i)d_i(m+1)d_i(m)
+W(m,i)d_i^{\,2}(m).
\end{align}
Hence the theorem says that  \eqref{form-1}
   is positive.
  If
  $\Delta_2(m,i)<0$, it
   is obvious that  \eqref{form-1} is positive
     since
   $U(m,i)>0$ for $1\leq i \leq m-1$.
   We now assume that $\Delta_2(m,i)\geq 0$.

   Recall  the
   five cases for the range of $i$ as given before.  Case 1:
 $1\leq i < \left(\frac{m^2}
 {2}\right)^{1/3}-m^{1/3} $.  By Theorem
  \ref{thm-4-6},
  we see that \eqref{form-1}  is
positive. Case 2: $\left(\frac{m^2}
 {2}\right)^{1/3}-m^{1/3}\leq i
 \leq \left(\frac{m^2}
 {2}\right)^{1/3}$.  Note that in this case, by Lemma \ref{Lem41},
we have $\Delta_2(m,i)<0$, which belongs to the case that
 we have already considered before.
  Case 3:  $\left(\frac{m^2}{2}\right)^{1/3}
 < i < m^{2/3}$. It follows from  Theorem
   \ref{Lemma-4-7} that \eqref{form-1}
   is positive. Case 4: $m^{2/3}\leq i
    \leq m-4$. The lower bound given in
    Theorem \ref{thm-4-4}  ensures that
     \eqref{form-1} is positive.
    It remains to consider the case when $i=m-3,\; m-2,\; m-1$.
  Here we only verify the statement for $i=m-3$.  The
   other two cases can be justified analogously. By \eqref{Defi}, we see that
\begin{align*}
&U(m,m-3)d_{m-3}^{\,2}(m+1)
 +V(m,m-3)d_{m-3}(m+1)d_{m-3}(m)
\\[6pt]
&\qquad +W(m,m-3)d_{m-3}^{\,2}(m) =\frac{(m+1)^2(m-2)g(m)}
{9216(2m+1)^2(2m-1)^2(2m-3)^2} 2^{-2m} {2m+2 \choose m+1}^2,
\end{align*}
where $g(m)$ is given by
\begin{align*}
g(m)=&\;2048m^{12}-10240m^{11}+16512m^{10}
-3456m^9-35232m^8+99120m^7+44488m^6\\[6pt]
&-375620m^5 +431652m^4 -182601m^3+7362m^2
+13797m-2430,
\end{align*}
which is positive for $m\geq 273$.
This completes the proof. \qed

To conclude this paper, we show that the $2$-log-concavity
of the Boros-Moll polynomials implies
 the log-concavity of
the sequence $\{ i
   (i+1)(d_i^2(m)-d_{i-1}(m)d_{i+1}(m)) \}_{1\leq i \leq m}$, as
   stated in Theorem 1.7.

Clearly, for
$i\geq 2$, we have
 \begin{align}\label{I-I}
\frac{i(i+1)}{(i-1)(i+2)}>1.
 \end{align}
By Theorem \ref{Theo-1} and the inequality \eqref{I-I}, we obtain
that  for  $2\leq i \leq m-1$,
 \[
\frac{d_{i-1}^2(m)-d_{i-2}(m)d_{i}(m)} {
d_{i}^2(m)-d_{i-1}(m)d_{i+1}(m)} < \frac{i(i+1)}{(i-1)(i+2)}
\frac{d_{i}^2(m)-d_{i-1}(m)d_{i+1}(m)} {
d_{i+1}^2(m)-d_{i}(m)d_{i+2}(m)}.
 \]
Replacing $i$ by $i+1$, we find that for $1\leq i \leq m-2$,
 \[
\frac{d_{i}^2(m)-d_{i-1}(m)d_{i+1}(m)} {
d_{i+1}^2(m)-d_{i}(m)d_{i+2}(m)} < \frac{(i+1) (i+2)} {i(i+3)}
\frac{\left(d_{i+1}^2(m)-d_{i}(m)d_{i+2}(m)\right)} {
\left(d_{i+2}^2(m)-d_{i+1}(m)d_{i+3}(m)\right)},
 \]
which can be written  as
 \[
\frac{i(i+1)\left(d_{i}^2(m) -d_{i-1}(m)d_{i+1}(m)\right)} {
(i+1)(i+2)\left(d_{i+1}^2(m) -d_{i}(m)d_{i+2}(m)\right) } <
\frac{(i+1)(i+2)\left(d_{i+1}^2(m) -d_{i}(m)d_{i+2}(m)\right)}
{(i+2)(i+3)\left( d_{i+2}^2(m)
 -d_{i+1}(m)d_{i+3}(m)\right)}.
 \]
This means that the sequence  $\{ i
   (i+1)(d_i^2(m)-d_{i-1}(m)d_{i+1}(m)) \}_{1\leq i\leq m}$
   is log-concave.

\vspace{0.5cm}
 \noindent{\bf Acknowledgments.}
 This work was supported by
  the 973
Project, the PCSIRT Project of the Ministry of Education, the
Ministry of Science and Technology, and the National Science
Foundation of China.

\end{document}